\documentclass{article}
\usepackage{graphicx} 
\usepackage{amsmath, amsxtra, verbatim}
\usepackage{tabularx}

\usepackage[affil-it]{authblk}
\usepackage[mathletters]{ucs}
\usepackage[utf8x]{inputenc}

\usepackage{mathtools}
\usepackage[english]{babel}
\usepackage{amsthm}
\usepackage{amssymb}
\newtheorem{theorem}{Theorem}[section]
\newtheorem{prop}{Proposition}
\newtheorem{lemma}[theorem]{Lemma}
\newtheorem{definition}{Definition}[section]
\usepackage{graphicx}

\usepackage[cmtip,all]{xy}

\usepackage[colorlinks, linkcolor=blue]{hyperref}

\usepackage[shortlabels]{enumitem}

\def\RR{\mathbb{R}}
\def\PP{\mathbb{P}}

\def\NN{\mathbb{N}}
\def\EE{\mathbb{E}}

\title{Preventing Finite-Time Blowup in a Constrainted Potential for Reaction-Diffusion Equations}
\author{John Ivanhoe and Michael Salins}
\affil{Boston University}
\date{\today}

\begin{document}

\maketitle

\section{Introduction}

We are interested in studying how dissipative forcing prevents explosion of mild solutions to the stochastic reaction-diffusion equation (SRDE) defined on an open bounded domain $D \subset \RR^d$ with appropriately-smooth boundary. Our equation of interest is
\begin{equation}
    \label{eq:SRDE} \begin{cases}
    \frac{\partial u}{\partial t}(t,x) = \mathcal{A}u(t,x) + f(u(t,x)) + \sigma(u(t,x))\dot{W}(t,x) \\
    u(t,x) = 0, \;\; x \in \partial D, \hspace{0.25in} u(0,x) = u_0(x)
    \end{cases}
\end{equation}

Here, $\mathcal{A}$ is a second-order elliptic differential linear operator and $\dot{W}$ is a Gaussian noise. The function $f(u(t,x))$ represents a state-dependent external forcing and the multiplicative noise $\sigma(u(t,x))\dot{W}$ represents state-dependent stochastic forcing. We are particularly interested when $f$ models a constrained potential force, where $\lim_{|x| \to 1} f(x) \text{sign}(x) = -\infty$. We also allow $\sigma$ to be unbounded near $|x| \to 1$ and identify a sufficient condition on the relative growth rates of $f$ and $\sigma$ as $|x| \to 1$ that guarantees the existence of global solutions. \par

Specifically, we assume there are constants $c_0 \in (0,1)$, $C,K,\beta > 0$, and $\gamma \ge 0$ such that for $c_0 < |w| < 1$
\begin{equation}
\label{eq:growthrates} \begin{cases} |\sigma(w)| \le C(1-|w|)^{-\gamma} \\ f(w) \; \text{sign}(w) \le -K(1-|w|)^{-\beta} \end{cases}
\end{equation}

With $\gamma = 0$, we allow for additive noise. Analogous to \cite{salins2021global}, the stochastic forcing pushes solutions toward the endpoints with growing intensity, counteracting the dissipative force of $f$. 

This set-up is motivated by ecological models. In population growth models, a common assumption is letting $f$ and $\sigma$ obey logistic growth. One such example is the FKPP equation on the circle \cite{ShigaKFPP}:
\[
\partial_t u = \frac{\alpha}{2} \Delta u + \beta u(1-u) + \sqrt{\gamma u(1-u)} \dot{W}
\]

This FKPP model has been studied by various authors \cite{fan2024quasistationary}, \cite{DOERING2003243}, \cite{HobsonTribeFKPP}, \cite{Shiga_1994fkpp}, \cite{Mueller2009fkpp}, \cite{mueller2009effectfkpp}, with applications recorded in the survey \cite{Panja_2004}. The terms $u(1-u)$ assume logistic growth, where a population approaching its carrying capacity sees the forcing terms vanish. However, extreme downward forcing is also reasonable for a population nearing its upper limit (for example, if there are more elk than what a forest can sustain, many elk should die off). The stochastic forcing should also be large, reflecting the large, yet random, amount of offspring expected from a large population, as well as modeling uncertainty as we near physical limits. Together, these ideas form our assumptions in \eqref{eq:growthrates}. \par

With these assumptions, the SRDE \eqref{eq:SRDE} is not well-defined whenever $|u(t,x)|$ $\ge 1$, motivating our interest in preventing finite-time blowup, which we intuitively understand as any point $(t,x) \in (0, \infty) \times D$ such that $|u(t,x)| = 1$, assuming $||u(0,x)||_{L^\infty(D)} < 1$. We will show that, among other assumptions, the inequality
\begin{equation}
\label{eq:KeyAssumption} \gamma + 1 < \frac{(1-\eta)(\beta+1)}{2}
\end{equation}
prevents finite-time blowup, where $\gamma, \beta$ are our constants from \eqref{eq:growthrates}. The constant $\eta \in [0,1)$ was introduced by Cerrai (2003) and describes the balance between the roughness of the noise and smoothing of the elliptic operator \cite{Cerrai2003}. A larger $\eta$ corresponds to rougher solutions, with $\frac{1-\eta}{2}$ being the temporal Hölder continuity of the paths. This is reflected in inequality \eqref{eq:KeyAssumption}, where rougher solutions require larger $\beta$ values to prevent blow-up for a fixed $\gamma$. For an explicit definition of $\eta$, see \text{\bf{Assumption 3}}. \par

This result is analogous to the results in \cite{salins2021global}, where solutions to \eqref{eq:SRDE} are allowed to grow arbitrarily large under superlinear forcing. There, the key inequality preventing finite-time blowup is $\gamma - 1 < \frac{(1-\eta)(\beta-1)}{2}$. The two results are similar, differing largely due to the different exponents one gets from integrating $u^{-\beta}$ compared to $u^\beta$. However, our inequality presents a strange question in the case of additive noise. When $\gamma = 0$, our main result shows that finite-time blow-up is prevented for a strictly positive $\beta > \frac{2}{1-\eta} - 1 > 0$. It is unknown if this $\beta$ is critical in general, though criticality is known in the case of space-time white noise in one spatial dimension \cite{Mueller1999}. Criticality would imply that there exist exponents $\beta > 0$ such that finite-time blow up is possible for additive noise. In Mueller's investigation of $\eqref{eq:SRDE}$ with $\mathcal{A} = \Delta, f(u) = u^{-\beta}$, additive noise, and spatial dimension $d = 1$, it was shown that solutions couldn't decrease to $0$ in finite time (what we refer to as blow-up) when $\beta > 3$ \cite{Mueller1998LongtimeEF}. In Mueller's case, $\eta = \frac{1}{2}$ and $\gamma = 0$, making his result equivalent to \eqref{eq:KeyAssumption}. Mueller and Pardoux later showed that $\beta = 3$ was critical: for any $0 < \beta < 3$, there is a positive probability for the solution to hit $0$ in finite time. \par

Another example where solutions can blow-up in finite time is the $1$ dimensional SDE
\[
dX(t) = X(t)^{-\beta} dt + X(t)^{-\gamma} dB_t
\]

It can be shown (see \cite{cerrai2008khasminskii}) that finite-time blowup is prevented so long as
\[
\gamma+1 < \frac{\beta+1}{2},
\]

For Brownian motion, $\eta = 0$, so this condition aligns with our inequality \eqref{eq:KeyAssumption}. It is currently unknown if \eqref{eq:KeyAssumption} is the optimal condition for superlinear forcing terms and colored noise, and this problem is left for future work.

In the study of $\eqref{eq:SRDE}$, it is known that stochastic forcing can result in solutions growing arbitrarily large with positive probability. Mueller and collaborators \cite{Mueller1991longtime, Mueller1998LongtimeEF, mueller2000critical,Mueller1993blowup} investigated the case where $\mathcal{A} = \Delta, f(u) = 0$, and the spatial dimension $d = 1$. It was found that solutions can exhibit finite-time blowup if $|\sigma(u)| > c|u|^\gamma$ for some $c > 0$ and $\gamma > \frac{3}{2}$. It was also observed that explosion is prevented if $|\sigma(u)| \le C(1+|u|^\gamma)$ for some $C > 0$ and $\gamma < \frac{3}{2}$. See \cite{Pavel2018, PaoLiu2009,PaoLiu2011,foondun2019nonexistence} for investigations of explosion by other authors. \par
Early work showed that \eqref{eq:SRDE} has unique solutions if $f$ and $\sigma$ are globally Lipschitz continuous with at most linear growth \cite{daprato_zabczyk_2014, dalang2008minicourse, Dalang1999, Peter1992, liu2015stochastic, Peszat2000, MartaMonica2000, Sowers1992, Walsh1986AnIT}. However, in the case of superlinearly growing $\sigma$, a linearly growing $f$ is no longer sufficient to prevent explosive growth. Cerrai demonstrated the existence of global solutions for when $\sigma$ is locally Lipschitz continuous with linear growth and $f$ is strongly dissipative with polynomial growth \cite{Cerrai2003}. Assuming polynomial growth for $f$ is common \cite{Zdzislaw1999, Cerrai2011,Iwata1987,Manthey1999}, but Da Prato and Röckner \cite{DaPratoRockner2002} and Marinelli and Röckner \cite{marinelli2010uniqueness} proved that this polynomial growth restriction can be relaxed, such that existence and uniqueness of solutions is implied by a monotonicity condition on $f$. See also \cite{Cerrai2013, gordina2020ornsteinuhlenbeck, marinelli2018wellposedness, marinelli2010uniqueness}. Recently, it was shown that superlinearly-growing dissipative forcing $f$ can prevent blow-up of solutions in the case of superlinearly-growing $\sigma$ \cite{salins2021global}. All of these references consider the case where $f: \RR \to \RR$ and $u(t,x)$ can take any value in $\RR$.  \par

Since we assume our constraining force $f$ is of the form $f(w) = -\text{sign}(w)(1-|w|)^{-\beta}$ for some $\beta > 0$, we naturally assume $||u_0||_{L^\infty(D)} < 1$. In the absence of stochastic forcing, $f$ restricts $u(t,x)$ to remain in $(-1,1)$. Since equation \eqref{eq:SRDE} breaks down whenever $|u(t,x)| = 1$, this motivates the following definition of explosion: for each $n \in \NN$, we define the stopping times
\begin{equation}
\label{eq:LocalStoppingTimes} T_n = \inf \left\{t > 0 : ||u(t)||_{L^\infty(D)} > 1 - \frac{1}{3^n} \right\}
\end{equation}

Under $\text{\bf{Assumption 2}}$, which can be found in Section 2, below, the realization $A$ of $\mathcal{A}$ in $L^2(D)$ generates a $C_0$-semigroup $S(t)$. \par

\begin{definition}
A $C_0(\overline{D})$-valued process $u(t)$ is a local mild solution to $\eqref{eq:SRDE}$ if
\begin{equation}
\label{eq:MildSol} u(t) = S(t)u_0 + \int_0^t S(t-s) f(u(s)) ds + \int_0^t S(t-s) \sigma(u(s)) dW(s)
\end{equation}

for all $t \in [0,T_n]$ for any $n$. \par
\end{definition}

\begin{definition}
A mild solution $u$ is global if $u(t)$ solves $\eqref{eq:MildSol}$ for all $t > 0$ with probability $1$.
\end{definition}

Local mild solutions $u$ are continuous by definition, so $T_n \le T_{n+1}$ for all $n \in \NN$. Therefore, a local mild solution $u$ solves $\eqref{eq:MildSol}$ on the time interval $[0, \sup_n T_n)$ with probability $1$. Thus, a local mild solution $u$ is a global solution if and only if $\sup_n T_n = \infty$ with probability 1. We define the explosion time as $\sup_n T_n$, as mild solutions to \eqref{eq:SRDE} do not make sense after $\sup_n T_n$, and can now state our key result: \par

\begin{theorem} \label{thm:GSol} 
Let $\frac{1-\eta}{2}$ be the time Hölder continuity of our solution to \eqref{eq:SRDE}, with $\eta$ defined in Section 2. If Assumptions 1-4 are satisfied and $\gamma + 1 < \frac{(1-\eta)(\beta+1)}{2}$, then any local mild solution is a global solution.
\end{theorem}

In Section 2, we describe our notation, main assumptions, and main result. In Section 3, we highlight some crucial estimates and prove Lemma $\eqref{lem:Limgrowthlemma}$, which is central to limiting the growth of our solutions. In Section 4, we prove our main theorem.

\section{Assumptions and Main Result}

Let $D \subset \RR^d$ denote an open, bounded domain. For $p \in [1, \infty)$, define $L^p(D)$ to be the Banach space of functions $v: D \to \RR$ such that the norm
\[
|v|_{L^p(D)} = \left( \int_D |v(x)|^p dx \right)^\frac{1}{p}
\]
is finite. When $p = \infty$, the $L^\infty(D)$ norm is
\[
|v|_{L^\infty(D)} = \sup_{x \in D} |v(x)|
\]

Define $C_0(\overline{D})$ to be the subset of $L^\infty(D)$ of continuous functions $v: \overline{D} \to \RR$ such that $v(x) = 0$ for $x \in \partial D$. Define $C_0([0,T] \times \overline{D})$ to be the set of continuous functions $v: [0,T] \times \overline{D} \to \RR$ such that $v(t,x) = 0$ for $x \in \partial D$, endowed with the supremum norm
\[
|v|_{C_0([0,T] \times \overline{D})} = \sup_{t \in [0,T]} \sup_{x \in \overline{D}} |v(t,x)|
\]

We make the following assumptions about the differential operator $\mathcal{A}$, the noise $\dot{W}$, and the forcing terms $f$ and $\sigma$. 

\vspace{0.1in}

\text{\bf{Assumption 1}} The initial data $u_0 \in C_0(\overline{D})$ and $||u_0||_{\infty} < 1$

\vspace{0.1in}

\text{\bf{Assumption 2}} $\mathcal{A}$ is a self-adjoint second-order elliptic differential operator
\[
\mathcal{A} \phi(x) = \sum_{i=1}^d \sum_{j=1}^d \frac{\partial}{\partial x_j}\left( a_{ij}(x) \frac{\partial \phi}{\partial x_j} (x) \right) 
\]
where $a_{ij}$ are continuously differentiable on $\overline{D}$, symmetric, and elliptic. Let $A$ be the realization of $\mathcal{A}$ in $L^2(D)$ with Dirichlet boundary conditions. Then there exists a sequence of eigenvalues $0 \le \alpha_1 \le \alpha_2 \le \ldots $ and eigenfunctions $e_k \in L^2(D) \cap C_0(D)$ such that
\[
A e_k = -\alpha_k e_k \hspace{0.25in} |e_k|_{L^2(D)} = 1
\]

We believe our results hold for more general second-order elliptic differential operators. However, for convenience, we assume self-adjoint operators. Using \text{\bf{Assumption 2}} and estimates similar to those in Cerrai 2003 \cite{Cerrai2003}, we have an elliptic contraction semigroup $S(t)$ generated by this operator $A$, where $S(t)$ acts on functions $u(t,x)$ through the spatial integral
\begin{equation}
\label{eq:semigroup} (S(t)u(s))(x) = \int_D K(t,x,y)u(s,y) dy
\end{equation}

where $K(t,x,y)$ is the kernel of this contraction semigroup described by
\begin{equation}
\label{eq:semigroupkernel} K(t,x,y) = \sum_{k=1}^\infty e^{-\alpha_k t} e_k(x) e_k(y)
\end{equation}

\vspace{0.1in}

\text{\bf{Assumption 3}} (See \cite{Cerrai2003}, \cite{cerrai2008khasminskii}) There exists a sequence of numbers $\lambda_j \ge 0$ and a sequence of i.i.d one-dimensional Brownian motions $\{B_j(t)\}_{j \in \NN}$ such that, formally,
\[
\dot{W}(t,x) = \sum_{k=1}^\infty \lambda_j e_j(x) dB_j(t)
\]

Furthermore, there exist exponents $\theta > 0$ and $\rho \in [2, \infty)$ such that
\begin{align*}
    &\begin{cases}
    \left( \sum_{j=1}^\infty \lambda_j^\rho |e_j|_{L^\infty(D)}^2 \right)^\frac{2}{\rho} < \infty & \text{if $\rho \in [2,\infty)$} \\
    \sup_j \lambda_j < \infty & \text{if $\rho = \infty$}
    \end{cases} \\
    &\sum_{k=1}^\infty \alpha_k^{-\theta} |e_k|_{L^\infty(D)}^2 < \infty
\end{align*}
and
\[
\eta \coloneqq \frac{\theta(\rho-2)}{\rho} < 1
\]

This constant $\eta < 1$ represents how the noise $\dot{W}$ and elliptic operator $\mathcal{A}$ relate to each other, where $\eta < 1$ represents the smoothing of our elliptic operator being strong enough to counteract the irregularity of the noise so that solutions are function-valued. A trace-class noise $\eta = 0$ corresponds to $\rho = 2$. For space-time white noise on a one-dimensional spatial interval, it is known that $\alpha_k \sim k^2$, so we can take $\rho = \infty$ and $\eta = \theta$ can be any number larger than $\frac{1}{2}$.

\vspace{0.1in}

\text{\bf{Assumption 4}} The functions $f:(-1,1) \to \RR$ and $\sigma:(-1,1) \to \RR$ are locally Lipschitz continuous functions. There exist powers $\beta > 0$, $\gamma \ge 0$ and constants $c_0, C, K > 0$ such that for $|u| \in (c_0, 1)$
\begin{align*}
|\sigma(u)| &\le C(1-|u|)^{-\gamma} \\
f(u) \text{sign(u)} &\le -K(1-|u|)^{-\beta}
\end{align*}

Our main goal is in showing that local mild solutions can be global solutions, but this is only meaningful if local mild solutions exist. To show this, for each $n \in \NN$, define functions $f_n, \sigma_n: \RR \to \RR$ by
\begin{align*}
f_n(w) &= \begin{cases} f(-1 + \frac{1}{3^n}) & \text{if $w < -1 + \frac{1}{3^n}$} \\ f(w) & \text{if $|w| < 1 - \frac{1}{3^n}$} \\ f(1 - \frac{1}{3^n}) & \text{if $w > 1 - \frac{1}{3^n}$} \end{cases} \\
\sigma_n(w) &= \begin{cases} \sigma(-1 + \frac{1}{3^n}) & \text{if $w < -1 + \frac{1}{3^n}$} \\ \sigma(w) & \text{if $|w| < 1 - \frac{1}{3^n}$} \\ \sigma(1 - \frac{1}{3^n}) & \text{if $w > 1 - \frac{1}{3^n}$} \end{cases}
\end{align*}

These functions $f_n, \sigma_n$ are globally Lipschitz continuous functions. Thus, the SRDE
\begin{equation}
    \label{eq:CutoffSRDE} \begin{cases}
    \frac{\partial u_n}{\partial t}(t,x) = \mathcal{A}u_n(t,x) + f_n(u_n(t,x)) + \sigma_n(u_n(t,x))\dot{W}(t,x) \\
    u_n(t,x) = 0, \;\; x \in \partial D, \hspace{0.25in} u_n(0,x) = u_0(x)
    \end{cases}
\end{equation}

has a unique solution by the usual existence and uniqueness theorems \cite{daprato_zabczyk_2014}. Since $u_n$ solves \eqref{eq:CutoffSRDE}, it solves \eqref{eq:SRDE} up to $T_n$ due to our definitions of $f_n, \sigma_n$. Since $u_n$ uniquely solves \eqref{eq:CutoffSRDE}, it is a unique local solution. Furthermore, by this uniqueness and our construction of $f_n, \sigma_n$, the solutions $u_n$ are consistent, where $u_n = u_m$ until $T_n$ whenever $n < m$.

\vspace{0.1in}

Let us recall the definition of a mild solution:
\begin{definition}
A $C_0(\overline{D})$-valued process $u(t)$ is a local mild solution to $\eqref{eq:SRDE}$ if
\[
 u(t) = S(t)u_0 + \int_0^t S(t-s) f(u(s)) ds + \int_0^t S(t-s) \sigma(u(s)) dW(s)
\]

for all $t \in [0,T_n]$ for any $n$, where

\[
T_n = \inf \left\{t > 0 : ||u(t)||_{L^\infty(D)} > 1 - \frac{1}{3^n} \right\}
\]

\end{definition}

Under our assumptions, these stopping times are well-defined and we can state our final result

\begin{theorem} \label{thm:UGSol} 
Under Assumptions 1-4, there exists a unique, global solution to $\eqref{eq:SRDE}$.
\end{theorem}

\section{Estimates}

Using the factorization method of Da Prato and Zabczyk [Chapter 5.3.1], see also Cerrai (2003), the stochastic integral 
\begin{equation}
\label{eq:Zdef} Z(t) = \int_0^t S(t-s) \sigma(u(s)) dW(s),
\end{equation}

where $S(t)$ is the elliptic semigroup defined in \eqref{eq:semigroup}, can be written as 
\begin{equation}
\label{eq:FactorizationPrato} Z(t) = \frac{\sin(\pi \alpha)}{\pi} \int_0^t (t-s)^{\alpha - 1} S(t-s)  Z_\alpha(y) dy
\end{equation}

for $\alpha \in (0,1)$ and
\begin{equation}
\label{eq:ZalphaEq} Z_\alpha(s) = \int_0^t (t-s)^{-\alpha} S(t-s) \sigma(u(s)) dW(s)
\end{equation}

It is then possible to write $Z(t \wedge \tau)$, where $\tau$ is any stopping time with respect to the natural filtration of $W(t)$, as
\begin{align}
Z(t \wedge \tau)& = \frac{\sin(\pi \alpha)}{\pi} \int_0^{t \wedge \tau} (t \wedge \tau-s)^{\alpha - 1} S(t \wedge \tau-s)  Z_\alpha(s) \nonumber ds \\
\label{eq:StoppingTimeFactorization} &=\frac{\sin(\pi \alpha)}{\pi} \int_0^{t \wedge \tau} (t \wedge \tau-s)^{\alpha - 1} S(t \wedge \tau-s)  \tilde{Z}_\alpha(s) ds
\end{align}

where
\begin{equation}
\label{eq:IndicatorZalphaEq} \tilde{Z}_\alpha(t) = \int_0^t (t-s)^{-\alpha} S(t-s) \sigma(u(s)) 1_{\{s \le \tau\}} dW(s)
\end{equation}

We will make use of the following propositions

\begin{prop}
Let $\alpha \in \left(0, \frac{1-\eta}{2} \right)$ and $p \ge 2$. For any $t > 0$,
\begin{equation}
    \EE|\tilde{Z}_\alpha(t)|_{L^p(D)}^p \le C_{\alpha,p} \EE\left( \int_0^t (t-s)^{-\eta-2\alpha} |\sigma(u(s))|_{L^\infty(D)}^2 1_{ \{s \le \tau\} } ds \right)^{\frac{p}{2}}
\end{equation}
\end{prop}

Proposition 1 and its proof can be found in \cite{salins2021global} under Proposition (3.1), with its proof in that paper's Appendix. 

\begin{prop} For $p,\alpha, \beta$ such that $p(\alpha - \frac{1}{\beta+1}) - \frac{d}{2} - 1 > 0$, we have
\begin{equation}
\label{eq:ModifiedZMomentBound} \EE \left( \sup_{t \in [0,T]} \sup_{x \in D} \frac{|Z(t,x)|^p}{t^{\frac{p}{\beta+1}}}\right) \le C_{\alpha, p, d} T^{p(\alpha - \frac{1}{\beta+1}) - \frac{d}{2} - 1} \int_0^T \EE |\tilde{Z}_\alpha(s)|_{L^p(D)}^p ds
\end{equation}
\end{prop}

\emph{Proof of Proposition 2} Consider the factorized form of our stochastic convolution,
\begin{equation}
\label{eq:ZKernelForm} Z(t,x) = \frac{\sin(\pi \alpha)}{\pi} \int_0^t \int_D (t-s)^{\alpha - 1} K(t-s, x, y) \tilde{Z}_\alpha(s,y) dy ds
\end{equation}

where $K(t,x,y)$ is the kernel of the semigroup $S(t)$ defined in \eqref{eq:semigroupkernel}. Applying Holder's Inequality with $p$ and $\frac{p}{p-1}$ in both space and time implies
\begin{align}
|Z(t,x)| &\le \frac{\sin(\pi \alpha)}{\pi} \left( \int_0^t \int_D (t-s)^{ \frac{(\alpha - 1)p}{p-1}} |K(t-s, x, y)|^{\frac{p}{p-1}} dy ds \right)^{\frac{p-1}{p}} \nonumber \\
&\hspace{0.5in} \cdot \left( \int_0^t \int_D |\tilde{Z}_\alpha(s,y)|^p dy ds  \right)^{\frac{1}{p}} \nonumber \\
\label{eq:ZHolderIneq} &= \frac{\sin(\pi \alpha)}{\pi} \left( \int_0^t (t-s)^{ \frac{(\alpha - 1)p}{p-1}} ||K(t-s, x, \cdot )||_{L^{\frac{p}{p-1}}(D)}^{\frac{p}{p-1}} ds \right)^{\frac{p-1}{p}}  \\ 
&\hspace{0.5in} \cdot \left( \int_0^t \int_D |\tilde{Z}_\alpha(s,y)|^p dy ds\right)^{\frac{1}{p}} \nonumber
\end{align}

Using Holder's Inequality with the semigroup norm implies
\begin{align}
\label{eq:SemigroupHolderIneq} &\le \frac{\sin(\pi \alpha)}{\pi} \left( \int_0^t (t-s)^{ \frac{(\alpha - 1)p}{p-1}} ||K(t-s, x, \cdot )||_{L^1(D)} ||K(t-s,x,\cdot)||_{L^\infty(D)}^{\frac{1}{p-1}} ds \right)^{\frac{p-1}{p}} \\
&\hspace{0.5in} \cdot \left( \int_0^t \int_D |\tilde{Z}_\alpha(s,y)|^p dy ds\right)^{\frac{1}{p}} \nonumber
\end{align}

Properties of the semigroup kernel norm imply for any $t > 0$ and $x \in D$ that
\begin{align}
&\le \frac{\sin(\pi \alpha)}{\pi} \left( \int_0^t (t-s)^{ \frac{(\alpha - 1)p}{p-1}} \cdot 1 \cdot (t-s)^{-\frac{d}{2(p-1)}} ds \right)^{\frac{p-1}{p}} \nonumber \\
&\hspace{0.5in} \cdot \left( \int_0^t \int_D |\tilde{Z}_\alpha(s,y)|^p dy ds\right)^{\frac{1}{p}} \nonumber \\
&= \frac{\sin(\pi \alpha)}{\pi} \left( \int_0^t s^{ \frac{(\alpha - 1)p - \frac{d}{2} }{p-1}} ds \right)^{\frac{p-1}{p}} \cdot \left( \int_0^t \int_D |\tilde{Z}_\alpha(s,y)|^p dy ds\right)^{\frac{1}{p}} \nonumber \\
 &= C_{\alpha,p,d} t^{ \frac{p\alpha - \frac{d}{2} - 1}{p} } \left( \int_0^t \int_D |\tilde{Z}_\alpha(s,y)|^p dy ds\right)^{\frac{1}{p}} \nonumber \\
\label{eq:PolynomialZBound} |Z(t,x)|^p &\le C_{\alpha,p,d} t^{p\alpha - \frac{d}{2} - 1} \int_0^t |\tilde{Z}_\alpha(s)|_{L^p(D)}^p ds
\end{align}

Since $\eqref{eq:PolynomialZBound}$ holds for any $t > 0$ and $x \in D$, we can divide both sides by $t^{\frac{p}{\beta+1}}$ to see that,  for any $t > 0$,
\[
\frac{|Z(t,x)|^p}{t^{\frac{p}{\beta+1}}} \le C_{\alpha,p,d} t^{p(\alpha - \frac{1}{\beta+1}) - \frac{d}{2} - 1} \int_0^t |\tilde{Z}_\alpha(s)|_{L^p(D)}^p ds
\]

Since supremums preserve inequalities, we can take expectations to get
\begin{align}
\EE \left( \sup_{t \in [0,T]} \sup_{x \in D} \frac{|Z(t,x)|^p}{t^{\frac{p}{\beta+1}}}\right) &\le C_{\alpha, p, d} \sup_{t \in [0,T]} \left( t^{p(\alpha - \frac{1}{\beta+1}) - \frac{d}{2} - 1} \int_0^t \EE |\tilde{Z}_\alpha(s)|_{L^p(D)}^p ds \right)\\
&= C_{\alpha, p, d} T^{p(\alpha - \frac{1}{\beta+1}) - \frac{d}{2} - 1} \int_0^T \EE |\tilde{Z}_\alpha(s)|_{L^p(D)}^p ds
\end{align}

so long as we choose $\beta,p$ such that $p(\alpha - \frac{1}{\beta+1}) > 1 + \frac{d}{2}$, causing the right-hand side to be increasing in $t$. \qed

\subsection{Uniform Bounds}

To analyze the mild solution to the SPDE, we first demonstrate how the growth of the mild solution can be controlled by a deterministic function. Consider our mild solution to the integral equation
\begin{equation}
\label{eq:SPDESolution} u(t) = S(t)u_0 + \int_0^t S(t-s) f(u(s)) ds + Z(t)
\end{equation}

where we have suppressed the spatial variable for convenience. Here, $Z \in C_0([0,T] \times \overline{D})$ is the continuous stochastic convolution integral from $\eqref{eq:Zdef}$. For notational convenience, we also define 
\begin{align}
    \label{eq:disteq} e(t) \coloneqq 1 - |u(t)|_{L^\infty(D)},
\end{align}
so that $e(0) = 1 - |u_0|_{L^\infty(D)}$. Under this notation, explosion occurs if $e(t)$ ever hits 0.

\begin{lemma}\label{lem:Limgrowthlemma} 
Let $T > 0$ and assume that $Z,u \in C_0([0,T] \times D)$ solves $\eqref{eq:SPDESolution}$ such that, for all $t \in [0,T]$,
\begin{align} 
& \label{eq:smallZ} |Z(t)|_{L^\infty(D)} \le \frac{1}{3}e(t)  \\ 
&\label{eq:largeu} e(t) \in \left(0, \frac{1}{3^N} \right] 
\end{align}
Then for all $t \in [0,T]$,
\begin{equation}
 \label{eq:Limgrowtheq} e(t)\ge \frac{3}{4} \left(e(0)^{\beta + 1} +K \left(\frac{2}{5} \right)^\beta (1+\beta) t \right)^{\frac{1}{\beta+1}} > K_\beta t^{\frac{1}{\beta+1}}
\end{equation}

where $K_\beta > 0$ is a constant depending solely on $\beta$ and $N \in \NN$ is the minimal integer such that $\frac{2}{3^N} < 1 - c_0$, where $c_0$ is from Assumption 4. \par
\end{lemma}

\text{\emph{Proof}}.
If we define $v(t,x) = u(t,x) - Z(t,x)$ for $u,Z \in C_0([0,T] \times D)$, then $v$ is weakly differentiable and weakly solves the PDE 
\[
\frac{\partial v}{\partial t}(t,x) = \mathcal{A}v(t,x) + f(v(t,x)+Z(t,x))
\]
By a Yosida approximation (via Proposition 6.2.2 of \cite{Cerrai2001book} or Theorem 7.7 of \cite{daprato_zabczyk_2014}), we can assume that $v$ is a strong solution. By Proposition D.4 in the appendix of \cite{daprato_zabczyk_2014}, $t \mapsto |v(t)|_{L^\infty(D)}$ is left-differentiable such that
\begin{equation}
\label{eq:GrowthRateDerivative} \frac{d^-}{dt} |v(t)|_{L^\infty(D)} \le \mathcal{A}v(t,x_t)\text{sign}(v(t,x_t)) + f(v(t,x_t) + Z(t,x_t)) \text{sign}(v(t,x_t))
\end{equation}

where $x_t \in D$ is a maximizer satisfying 
\begin{equation}
\label{eq:maximizer} |v(t)|_{L^\infty(D)} = |v(t,x_t)| = v(t,x_t) \text{sign}(v(t,x_t))
\end{equation}

The ellipticity of $\mathcal{A}$ and convexity at a maximizer implies
\[
\mathcal{A}v(t,x_t)\text{sign}(v(t,x_t)) \le 0
\]

Thus we need only examine
\begin{align}
    \label{eq:vbound} \frac{d^-}{dt} |v(t)|_{L^\infty(D)} \le f(v(t,x_t) + Z(t,x_t)) \text{sign}(v(t,x_t))
\end{align}

Since $v(t) = u(t) - Z(t)$, using both versions of the triangle inequality combined with $\eqref{eq:smallZ}$ implies
\begin{align*}
|u(t)|_{L^\infty(D)} - |Z(t)|_{L^\infty(D)} \le &|v(t)|_{L^\infty(D)} \le |u(t)|_{L^\infty(D)} + |Z(t)|_{L^\infty(D)} \\
1 - |u(t)|_{L^\infty(D)} - |Z(t)|_{L^\infty(D)} \le 1 &- |v(t)|_{L^\infty(D)} \le 1 - |u(t)|_{L^\infty(D)} + |Z(t)|_{L^\infty(D)} \\
e(t) - |Z(t)|_{L^\infty(D)} \le 1 &- |v(t)|_{L^\infty(D)} \le e(t) + |Z(t)|_{L^\infty(D)}
\end{align*}

By Assumption $\eqref{eq:smallZ}$, this implies
\begin{align}
\label{eq:VUpperBound}1-|v(t)|_{L^\infty(D)} &\le \frac{4}{3}e(t)  \\
\label{eq:VLowerBound}1-|v(t)|_{L^\infty(D)} &\ge \frac{2}{3}e(t)
\end{align}

Using the fact that $|Z(t,x_t)| \le |Z(t)|_{L^\infty(D)}$ by definition of supremum, and $|v(t,x_t)| = |v(t)|_{L^\infty(D)}$, we similarly find that
\begin{align*}
|v(t)|_{L^\infty(D)} - |Z(t)|_{L^\infty(D)} \le &|v(t,x_t) + Z(t,x_t)| \le |v(t)|_{L^\infty(D)} + |Z(t)|_{L^\infty(D)} \\
1 - |v(t)|_{L^\infty(D)} - |Z(t)|_{L^\infty(D)} \le 1 &- |v(t,x_t) + Z(t,x_t)| \le 1 - |v(t)|_{L^\infty(D)} + |Z(t)|_{L^\infty(D)} \\
1 - |v(t)|_{L^\infty(D)} - \frac{1}{3}e(t) \le 1 &- |v(t,x_t) + Z(t,x_t)| \le 1 - |v(t)|_{L^\infty(D)} + \frac{1}{3}e(t)
\end{align*}

Using inequalities $\eqref{eq:VUpperBound}$ and $\eqref{eq:VLowerBound}$ on the last series of inequalities implies
\begin{align}
\label{eq:vzUpperBound} 1 - |v(t,x_t) + Z(t,x_t)| &\le \frac{5}{3} e(t) \\
\label{eq:vzLowerBound} 1-|v(t,x_t) + Z(t,x_t)| &\ge \frac{1}{3}e(t)
\end{align}

Since $0 < e(t) \le \frac{1}{3^N}$ by $\eqref{eq:largeu}$, inequality $\eqref{eq:VUpperBound}$ implies $1-|v(t)|_{L^\infty(D)} \le \frac{4}{3^{N+1}}$ and $\eqref{eq:VLowerBound}$ implies $1-|v(t)|_{L^\infty(D)} > 0$. Furthermore,
\[
|v(t)|_{L^\infty(D)} > c_0 \iff 1 -  |v(t)|_{L^\infty(D)} < 1-c_0
\]

Since $\frac{2}{3^N} < 1-c_0$, this implies
\[
1-|v(t)|_{L^\infty(D)} \le \frac{4}{3^{N+1}} = \frac{2}{3} \frac{2}{3^N} < \frac{2}{3} (1-c_0)
\]

Thus, $|v(t)|_{L^\infty(D)} \in (c_0, 1)$ for $t \in [0,T]$. A similar argument implies $|v(t,x_t) + z(t,x_t)| \in (c_0, 1)$ for $t \in [0,T]$. Furthermore, 
\[
\text{sign}(v(t,x_t)) = \text{sign}(v(t,x_t)+z(t,x_t)),
\]
as  $|z(t)|_{L^\infty(D)}$ is sufficiently small relative to $|u(t)|_{L^\infty(D)}$ and $|v(t)|_{L^\infty(D)}$. By $\eqref{eq:VLowerBound}$ and $\eqref{eq:vzUpperBound}$, we get
\begin{align}
\label{eq:vzMainBound} 1-|v(t,x_t) + z(t,x_t)| \le \frac{5}{2} (1 - |v(t)|_{L^\infty(D)})
\end{align}

Applying $\eqref{eq:vzMainBound}$ to $\eqref{eq:vbound}$ implies
\begin{align}
\label{eq:strongderivative} \frac{d^-}{dt} |v(t)|_{L^\infty(D)} \le -K \left(\frac{2}{5} \right)^\beta (1 - |v(t)|_{L^\infty(D)})^{-\beta} 
\end{align}

Letting $F(x) = \frac{1}{1+\beta}(1-x)^{\beta + 1}$ and $F'(x) = - (1-x)^\beta$, this implies
\begin{equation}
\label{eq:DerivativeFBound} \frac{d^-}{dt} F(|v(t)|_{L^\infty(D)}) \ge K \left(\frac{2}{5} \right)^\beta 
\end{equation}

Therefore, 
\begin{equation}
\label{eq:completedFBound}F(|v(t)|_{L^\infty(D)}) \ge F(|u_0|_{L^\infty(D)}) +K \left(\frac{2}{5} \right)^\beta t 
\end{equation}

and 
\begin{equation}
\label{eq:BoundOnV} |v(t)|_{L^\infty(D)} \le 1 - \left( e(0)^{\beta + 1} + K \left(\frac{2}{5} \right)^\beta (1+\beta) t \right)^{\frac{1}{\beta+1}}
\end{equation} 

with $\eqref{eq:BoundOnV}$ implying $ |v(t)|_{L^\infty(D)} \le |u_0|_{L^\infty(D)}$ for $t \in [0,T]$. Combining this result with $\eqref{eq:VUpperBound}$ yields our conclusion \eqref{eq:Limgrowtheq}:
\[
e(t) \ge \frac{3}{4} \left( e(0)^{\beta + 1} + K \left(\frac{2}{5} \right)^\beta (1+\beta) t \right)^{\frac{1}{\beta+1}}
\]
\qed \par

Inequality $\eqref{eq:Limgrowtheq}$ implies $|u(t)|_{L^\infty(D)} \in (0,1)$ for $t \in [0,T]$.  In addition, if we assume $1 - |u_0|_{L^\infty(D)} = \frac{1}{3^n}$ for some $n \in \NN$, then
\begin{align}
\label{eq:NoDecreaseResult} 1 - |u(t)|_{L^\infty(D)} \ge \frac{1}{4\cdot 3^{n-1}}  > \frac{1}{3^{n+1}}
\end{align} 
Thus, $e(t) > \frac{1}{3^{n+1}}$ for $t \in [0,T]$, meaning the distance between $1$ and $|u(t)|_{L^\infty(D)}$ can't decrease below $\frac{1}{3^{n+1}}$ so long as conditions $\eqref{eq:smallZ}$ and $\eqref{eq:largeu}$ hold. \par

\section{Proof of Theorem \ref{thm:GSol}/ Main Result}

To establish theorem \ref{thm:GSol}, we first define a sequence of stopping times. Recalling $e(t) = 1 - |u(t)|_{L^\infty(D)}$, we define
\begin{align}
& \tau_0 = \text{inf} \{t \ge 0: e(t) \le \frac{1}{3^n} \; \text{for some $n \in \{N,N+1,N+2,\ldots\}$ }\} \nonumber \\
\label{eq: LevelStoppingTimes}&\tau_{k+1} = \begin{cases} \text{inf}\{t \ge \tau_k : e(t) \le  \frac{1}{3^{N+1}} \} & \text{if $e(\tau_k) =  \frac{1}{3^N} $}\\ \text{inf}\{t \ge \tau_k : e(t) \le \frac{1}{3} e(\tau_k) \nonumber  \\ 
\hspace{0.5in} \text{or} \; e(t) \ge 3e(\tau_k)  \} & \text{if $e(\tau_k) \le \frac{1}{3^{N+1}} $}\end{cases} \\
\end{align}

Here, $N$ is the constant defined in Lemma \ref{lem:Limgrowthlemma}. While these $\tau_k$ aren't the same as the stopping times $T_k$ defined in \eqref{eq:LocalStoppingTimes}, it is clear that $\sup_k \tau_k < \infty$ if and only if $\sup_k T_k < \infty$.

\begin{lemma} \label{lem:BoundedProbabilityLemma} There exist constants $C > 0$ and $q > 1$, independent of $n,k,$ and $\epsilon > 0$, such that for any $\epsilon > 0$, any $k \in \NN$, and any $n \in \{N+1,N+2,\ldots\}$,
\[
\PP \left(e(\tau_{k+1}) = \frac{1}{3}e(\tau_k) \; \text{and} \; \tau_{k+1} - \tau_k < \epsilon \; \middle \vert \; e(\tau_k) = \frac{1}{3^n} \right) \le C \epsilon^q
\]
where $N$ is the constant defined in Lemma \ref{lem:Limgrowthlemma}.
\end{lemma}

\emph{Proof}
As discussed in $\bf{Lemma \; 4.1}$ from \cite{salins2021global}, we may assume $k = 0$, so $\tau_0 = 0$, without loss of generality by the Markov property. Therefore,
\begin{align*}
&\PP \left(  e(\tau_{k+1}) = \frac{1}{3}e(\tau_k) \; \text{and} \; \tau_{k+1} - \tau_k < \epsilon \; \middle \vert \; e(\tau_k) =\frac{1}{3^n} \right) \\
&= \PP \left(  e(\tau_1) = \frac{1}{3}e(0) \; \text{and} \; \tau_1 < \epsilon \; \middle \vert \; e(0) = \frac{1}{3^n} \right)
\end{align*}

Since $n \ge N+1$, by the definition of $\tau_1$, if $e(0) = \frac{1}{3^n}$, then $e(t) \in (\frac{1}{3^{n+1}}, \frac{1}{3^{n-1}})$ on $[0,\tau_1)$, equalling one of the endpoints at $t = \tau_1$. Since $\frac{1}{3^{n-1}} \le \frac{1}{3^N}$, assumption $\eqref{eq:largeu}$ is satisfied on the interval $t \in [0,\tau_1]$. \par
Let $Z(t)$ be the stochastic integral in \eqref{eq:SPDESolution}. We introduce the additional stopping time
\begin{equation}
\label{eq:smallZStoppingTime}\tilde{\tau} = \text{inf} \{t > 0 : |Z(t)|_{L^\infty(D)} > \frac{1}{3} e(t) \} 
\end{equation}

By definition, $|Z(t)|_{L^\infty(D)} \le \frac{1}{3} e(t)$ for $t \in [0,\tilde{\tau}]$. Lemma $\ref{lem:Limgrowthlemma}$ is satisfied up to $\tilde{\tau}$, therefore if $e(0) = \frac{1}{3^n}$, then $e(t) > \frac{1}{3^{n+1}}$ on $[0, \tilde{\tau} \wedge \tau_1]$. Thus, in the event $e(\tau_1) = \frac{1}{3} e(0)$, this implies that $\tilde{\tau} \le \tau_1$. If we were to assume $\tau_1 < \tilde{\tau}$ and $e(\tau_1) = \frac{1}{3} e(0)$, then $e(t)$ decreases to $\frac{1}{3^{n+1}}$ within the time interval $[0,\tilde{\tau}]$, which is impossible by Lemma $\ref{lem:Limgrowthlemma}$. \par

Consequently,
\begin{align}
&\PP \left(  e(\tau_1) = \frac{1}{3}e(0) \; \text{and} \; \tau_1 < \epsilon \; \middle \vert \; e(0) = \frac{1}{3^n} \nonumber \right) \\
&\le \PP \left( \tilde{\tau} \le \tau_1 < \epsilon  \; |  \; e(0) = \frac{1}{3^n} \nonumber \right) \\
&\le \PP \left( \sup_{t \in [0,\tilde{\tau} \wedge \tau_1 \wedge \epsilon]} \frac{|Z(t)|_{L^\infty(D)}}{e(t)} \ge \frac{1}{3} \; \middle \vert \; e(0) = \frac{1}{3^n} \nonumber \right) \\
\label{eq:preChebyshev} &\le \PP \left( \sup_{t \in [0,\tilde{\tau} \wedge \tau_1 \wedge \epsilon]} \frac{|Z(t)|_{L^\infty(D)}}{K_{\beta} t^{\frac{1}{\beta+1}}} \ge \frac{1}{3} \; \middle \vert \; e(0) = \frac{1}{3^n} \right),
\end{align} 

where $\eqref{eq:preChebyshev}$ uses the fact that $e(t) \ge K_{\beta} t^{\frac{1}{\beta+1}}$ on $[0,\tilde{\tau} \wedge \tau_1 \wedge \epsilon]$ by Lemma $\ref{lem:Limgrowthlemma}$. Chebyshev's Inequality yields
\begin{equation}
\label{eq:ChebyshevInequality}  \le \frac{\EE \left(\sup_{t \in [0,\tilde{\tau} \wedge \tau_1 \wedge \epsilon]} \frac{|Z(t)|_{L^\infty(D)}^p  }{t^{ \frac{p}{\beta+1}}} \right)}{K_{\beta}^p \frac{1}{3^p} }
\end{equation}

Using $\eqref{eq:ModifiedZMomentBound}$, by \text{\bf{Assumption 4}}, the above is bounded by

\begin{equation}
\label{eq: SmallProbability} \le K_{\beta,p} C_{\alpha,p,d} \epsilon^{p(\alpha - \frac{1}{\beta+1})-\frac{d}{2}-1} \int_0^T \EE |\tilde{Z}_\alpha(s)|_{L^p(D)}^p ds  
\end{equation}

From Lemma $\ref{lem:Limgrowthlemma}$, on the time interval $t \in [0,\epsilon \wedge \tilde{\tau} \wedge \tau_1]$
\begin{align}
\EE |\tilde{Z}_\alpha (t)|_{L^p(D)}^p &\le C \EE \left( \int_0^t (t-s)^{-2\alpha - \eta} |\sigma(u(s))|_{L^\infty(D)}^2 1_{\{s \le \tau\} } ds \right)^{\frac{p}{2}} \nonumber \\
&\le C \EE \left( \int_0^t (t-s)^{-2\alpha - \eta} (1 - |u(s)|_{L^\infty(D)})^{-2\gamma} 1_{\{s \le \tau\} } ds \right)^{\frac{p}{2}} \nonumber \\
\label{eq:ZTildeModifiedBound} &\le C_{\gamma, \beta}  \left( \int_0^t (t-s)^{-2\alpha - \eta} s^{-\frac{2\gamma}{\beta + 1}} ds \right)^{\frac{p}{2}} 
\end{align}

We require $\beta > 0$ to satisfy
\begin{equation}
\label{eq: BetaGammaCondition1} \eta + \frac{2\gamma}{\beta + 1 } < 1 \iff \gamma < \frac{(1-\eta)(\beta + 1)}{2},
\end{equation}

which is guaranteed by \text{\bf{Assumption 4}}. Choose $2\alpha = 1 - \eta - \frac{2\gamma}{\beta + 1}$, so the properties of the Beta Integral ensure that $\eqref{eq:ZTildeModifiedBound}$ equals a constant
\begin{equation}
\label{eq: BetaConstant} \left(\frac{\pi}{\sin(\pi(2\alpha + \eta))}\right)^{\frac{p}{2}} 
\end{equation}

Therefore, $\EE |\tilde{Z}_\alpha (t)|_{L^p(D)}^p \le C_{\alpha, \gamma, \eta, p, \beta} $, a constant independent of $\epsilon, n, k$. This result gives an upper bound for $\eqref{eq: SmallProbability}$:
\[
 \le C_{\alpha, p, d,\gamma, \eta} \epsilon^{p( \frac{1-\eta}{2} - \frac{\gamma+1}{\beta+1} ) - \frac{d}{2}} = \tilde{C} \epsilon^q
\]

by defining $q = p( \frac{1-\eta}{2} - \frac{\gamma+1}{\beta+1} ) - \frac{d}{2}$. \text{\bf{Assumption $4$}} guarantees that
\[
\frac{1-\eta}{2} - \frac{\gamma+1}{\beta+1} > 0
\]

Choose $p$ sufficiently large to ensure $q > 1$. Thus,

\begin{align}
&\PP \left(  e(\tau_{k+1}) = \frac{1}{3} e(\tau_k) \; \text{and} \; \tau_{k+1} - \tau_k < \epsilon \; \middle \vert \; e(\tau_k) = \frac{1}{3^n} \right) \leq C \epsilon^q
\end{align}

where $C> 0, q > 1$ independent of $\epsilon, k, n$ and $\epsilon$ independent of $k,n$, proving Lemma $\ref{lem:BoundedProbabilityLemma}$. \qed \par

Now we prove Theorem \ref{thm:GSol}.

\emph{Proof of Theorem 1.1} Using the stopping times $\tau_k$ defined above and Lemma $\eqref{lem:BoundedProbabilityLemma}$, it is known that there exists $C > 0$ and $q > 1$ such that for any $k \in \NN$ and small $\epsilon > 0$

\[
\PP \left(e(\tau_{k+1}) = \frac{1}{3}e(\tau_k) \; \text{and} \; \tau_{k+1} - \tau_k < \epsilon \right) \le C \epsilon^q
\]

If we pick $\epsilon = \frac{1}{k}$ for sufficiently-large $k$, then 
\[
\PP \left(e(\tau_{k+1}) = \frac{1}{3}e(\tau_k) \; \text{and} \; \tau_{k+1} - \tau_k < \frac{1}{k}\right) \le C \frac{1}{k^q}
\]

Since $q > 1$,
\[
\sum_{k=1}^\infty \PP \left(e(\tau_{k+1}) = \frac{1}{3}e(\tau_k) \; \text{and} \; \tau_{k+1} - \tau_k < \epsilon \right) < \infty
\]

Thus, by the Borel-Cantelli Lemma, with probability one there exists a random index $N_0(\omega) > 0$ such that for all $k \ge N_0(\omega)$, either
\begin{equation}
\label{eq:BorelCantelliResult} \tau_{k+1} - \tau_k \ge \frac{1}{k} \hspace{0.15in} \text{or}  \hspace{0.15in} e(\tau_{k+1}) = 3 e(\tau_k)
\end{equation}


It is possible that $\tau_K = \infty$ for some $K \in \NN$, which would imply $\tau_k = \infty$ for $k \ge K$ by our definition \eqref{eq: LevelStoppingTimes} of $\tau_k$.  This would imply that $e(\tau_{k+1}) = \frac{1}{3} e(\tau_k)$ finitely-many times, making it impossible for $||u(t)||_{L^\infty(D)}$ to hit $1$. Clearly, this implies mild solutions are global. If $e(\tau_{k+1}) = 3e(\tau_k)$ for all $k \ge N_0(\omega)$, then $||u(t)||_{L^\infty(D)}$ is distancing from $1$, again implying mild solutions are global. Thus, we are left with examining the case when $\tau_{k+1} - \tau_k \ge \frac{1}{k}$ for finite $\tau_k$ and $e(\tau_{k+1}) = \frac{1}{3} e(\tau_k)$ occurring infinitely-many times.

From the definition of our $\tau_k$, $\max_{k \ge N_0(\omega)} e(\tau_k)$ is attained. Thus, we may choose $N_1(\omega) > N_0(\omega)$ such that for $k \ge N_1(\omega)$, we have $e(\tau_k) \le e(\tau_{N_1(\omega)})$. \par
Now, for any $m \ge N_1(\omega) > N_0(\omega)$, we have
\begin{equation}
\label{eq:STDifferenceLowerBound} \sum_{k = N_1(\omega)}^m (\tau_{k+1} - \tau_k) \ge \sum_{k=N_1(\omega)}^m \frac{1}{k} 1_{\{e(\tau_{k+1}) = \frac{1}{3} e(\tau_k) \}}
\end{equation}

By the definition of our $N_1$, there must be more steps where $e
(\tau_k)$ decreases than steps where it increases, such that for $m \ge N_1(\omega)$:
\begin{equation}
\label{eq:LowerBoundSteps} U_m(\omega) = \sum_{k=N_1(\omega)}^{m-1} 1_{\{e(\tau_{k+1}) = \frac{1}{3} e(\tau_k) \}} \ge \frac{m-N_1(\omega)}{2}
\end{equation}

At this point, our argument follows the rest of \cite{salins2021global}, starting from equation (4.31). Using the summation by parts formula and \eqref{eq:LowerBoundSteps}, we find
\begin{align*}
\sum_{k = N_1(\omega)}^m (\tau_{k+1} - \tau_k) &\ge \sum_{k=N_1(\omega)}^m \frac{1}{k} 1_{\{e(\tau_{k+1}) = \frac{1}{3} e(\tau_k) \}} \\
&= \sum_{k=N_1(\omega)}^m \frac{1}{k} (U_{k+1}(\omega) - U_k(\omega)) \\
&= \frac{U_{m+1}(\omega)}{m} - \sum_{k = N_1(\omega)+1}^m U_k(\omega) \left(\frac{1}{k} - \frac{1}{k-1} \right) \\
&= \frac{U_{m+1}(\omega)}{m} + \sum_{k = N_1(\omega)+1}^m U_k(\omega) \left(\frac{1}{k(k-1)} \right) \\
&\ge \frac{m-N_1(\omega)}{2m} + \sum_{k = N_1(\omega)+1}^m \left(\frac{k-N_1(\omega)}{2k(k-1)} \right)
\end{align*}

As $m$ tends to infinity, the sum diverges. Therefore, with probability one,
\[
\sum_{k = N_1(\omega)}^\infty (\tau_{k+1} - \tau_k) = \infty
\]

and solutions cannot explode in finite time. \qed

\bibliographystyle{plain}
\bibliography{refs.bib}

\end{document}